\documentclass[reqno]{amsart}

\usepackage[all]{xy}
\usepackage{graphicx}

\usepackage{amssymb}
\usepackage{amsmath}
\usepackage{mathrsfs}
\usepackage{epsfig}
\usepackage{amscd}
\usepackage{graphicx, color}

\def\E{\ifmmode{\mathbb E}\else{$\mathbb E$}\fi} 
\def\N{\ifmmode{\mathbb N}\else{$\mathbb N$}\fi} 
\def\R{\ifmmode{\mathbb R}\else{$\mathbb R$}\fi} 
\def\Q{\ifmmode{\mathbb Q}\else{$\mathbb Q$}\fi} 
\def\C{\ifmmode{\mathbb C}\else{$\mathbb C$}\fi} 
\def\H{\ifmmode{\mathbb H}\else{$\mathbb H$}\fi} 
\def\Z{\ifmmode{\mathbb Z}\else{$\mathbb Z$}\fi} 
\def\P{\ifmmode{\mathbb P}\else{$\mathbb P$}\fi} 
\def\T{\ifmmode{\mathbb T}\else{$\mathbb T$}\fi} 
\def\SS{\ifmmode{\mathbb S}\else{$\mathbb S$}\fi} 
\def\DD{\ifmmode{\mathbb D}\else{$\mathbb D$}\fi} 

\renewcommand{\b}{\beta}

\newcommand{\e}{\varepsilon}

\renewcommand{\o}{\omega}

\newcommand{\del}{\partial}

\newcommand{\ben}{\begin{enumerate}}
\newcommand{\een}{\end{enumerate}}
\newcommand{\be}{\begin{equation}}
\newcommand{\ee}{\end{equation}}
\newcommand{\bea}{\begin{eqnarray}}
\newcommand{\eea}{\end{eqnarray}}
\newcommand{\beastar}{\begin{eqnarray*}}
\newcommand{\eeastar}{\end{eqnarray*}}
\newcommand{\bc}{\begin{center}}
\newcommand{\ec}{\end{center}}

\theoremstyle{theorem}
\newtheorem{thm}{Theorem}[section]
\newtheorem{cor}[thm]{Corollary}
\newtheorem{lem}[thm]{Lemma}
\newtheorem{prop}[thm]{Proposition}

\theoremstyle{definition}

\newtheorem{rem}[thm]{Remark}

\newtheorem*{thm*}{Theorem}

\numberwithin{equation}{section}

\def\R{{\mathbb R}}

\def\E{{\mathbb E}}
\def\Z{{\mathbb Z}}
\def\C{{\mathbb C}}
\def\R{{\mathbb R}}
\def\P{{\mathbb P}}

\def\N{{\mathbb N}}

\def\11{{\mathbb I}}

\def\delbar{{\overline \partial}}

\def\C{\mathbb{C}}
\def\Z{\mathbb{Z}}

\def\T{\mathbb{T}}

\def\Q{\mathbb{Q}}

\def\E{\ifmmode{\mathbb E}\else{$\mathbb E$}\fi} 
\def\N{\ifmmode{\mathbb N}\else{$\mathbb N$}\fi} 
\def\R{\ifmmode{\mathbb R}\else{$\mathbb R$}\fi} 
\def\Q{\ifmmode{\mathbb Q}\else{$\mathbb Q$}\fi} 
\def\C{\ifmmode{\mathbb C}\else{$\mathbb C$}\fi} 
\def\H{\ifmmode{\mathbb H}\else{$\mathbb H$}\fi} 
\def\Z{\ifmmode{\mathbb Z}\else{$\mathbb Z$}\fi} 
\def\P{\ifmmode{\mathbb P}\else{$\mathbb P$}\fi} 
\def\SS{\ifmmode{\mathbb S}\else{$\mathbb S$}\fi} 
\def\DD{\ifmmode{\mathbb D}\else{$\mathbb D$}\fi} 

\def\R{{\mathbb R}}

\def\E{{\mathbb E}}
\def\Z{{\mathbb Z}}
\def\C{{\mathbb C}}
\def\R{{\mathbb R}}

\def\N{{\mathbb N}}

\def\JJ{{\mathcal J}}

\def\delbar{{\overline \partial}}

\def\b{\beta}

\def\e{\varepsilon}

\def\o{\omega}  

  \def\S{\Sigma}

\def\CB{{\mathcal B}}

\def\CF{{\mathcal F}}

\def\CH{{\mathcal H}}

\def\CJ{{\mathcal J}}

\def\CM{{\mathcal M}}

\def\darr#1{\raise1.5ex\hbox{$\leftrightarrow$}
\mkern-16.5mu #1}

\def\roughly#1{\raise.3ex\hbox{$#1$\kern-.75em
\lower1ex\hbox{$\sim$}}}

\def\opname#1{\mathop{\kern0pt{\rm #1}}\nolimits}

\def\supp{\operatorname{supp}}

\def\Aut{\operatorname{Aut}}

\begin{document}
\quad \vskip1.375truein

\def\mq{\mathfrak{q}}
\def\mp{\mathfrak{p}}
\def\mH{\mathfrak{H}}
\def\mh{\mathfrak{h}}
\def\ma{\mathfrak{a}}
\def\ms{\mathfrak{s}}
\def\mm{\mathfrak{m}}
\def\mn{\mathfrak{n}}
\def\mz{\mathfrak{z}}
\def\mw{\mathfrak{w}}
\def\Hoch{{\tt Hoch}}
\def\mt{\mathfrak{t}}
\def\ml{\mathfrak{l}}
\def\mT{\mathfrak{T}}
\def\mL{\mathfrak{L}}
\def\mg{\mathfrak{g}}
\def\md{\mathfrak{d}}
\def\mr{\mathfrak{r}}

\title[Embedding property]
{Embedding property of $J$-holomorphic curves in
Calabi-Yau manifolds for generic $J$}
\author{Yong-Geun Oh}
\address{Department of Mathematics, University of Wisconsin-Madison, Madison, WI, 53706
\indent \& Korea Institute for Advanced Study, Seoul, Korea}
\email{oh@math.wisc.edu}
\thanks{The senior author is partially supported
by the NSF grant \#DMS 0503954}

\author{Ke Zhu}
\address{Department of Mathematics, University of Wisconsin-Madison, Madison, WI, 53706
\indent \& Department of Mathematics, The Chinese University of
Hong Kong, Hong Kong} \email{kzhu@math.cuhk.edu.hk}

\begin{abstract}
In this paper, we prove that for a generic choice of tame (or
compatible) almost complex structures $J$ on a symplectic manifold
$(M^{2n},\omega)$ with $n \geq 3$ and
with its first Chern class $c_1(M,\omega) = 0$, all somewhere injective $J$-holomorphic maps from
any closed smooth Riemann surface into $M$ are \emph{embedded}.
We derive this result as a consequence of the general optimal 1-jet evaluation
transversality result of $J$-holomorphic maps in general symplectic manifolds
that we also prove in this paper.
\end{abstract}

\date{Revision ; Jan 31, 2009}
\maketitle

\tableofcontents

\section{Introduction}
\label{sec:intro}

Let $(M,\omega)$ be a symplectic manifold of real dimension $2n$. We
denote by $J$ an almost complex structure tame to $\omega$ and by
$\CJ_\omega$ the set of tame almost complex structures. It is
a classical fact \cite{gromov}, \cite{mcduff} that for a generic choice of $J$, any
{\emph somewhere injective} $J$-holomorphic curve is a \emph{smooth} point
in the moduli space of $J$-holomorphic curves : A $J$-holomorphic
curve $u: \Sigma \to M$ is called somewhere injective if there is a
point $z \in \Sigma$ such that
$$
du(z) \neq 0 \quad \mbox{and } \, u^{-1}(u(z)) = \{z\}.
$$
This fact has been a fundamental point in the definition of
Gromov-Witten invariants and the counting problem of $J$-holomorphic
curves. Recent development in the Gromov-Witten theory unravels
necessity of finer structure theorem on the image of $J$-holomorphic
curves. In particular a conjectural mathematical definition of
Gopakuma-Vafa invariant of BPS-count is closely related to the
number of \emph{embedded} $J$-holomorphic curves in Calabi-Yau
three-folds for a generic choice of $J$ \cite{konts}, \cite{pandh}.

Now description of the main results of this paper is in order.

Let $\Sigma$ be a connected closed smooth surface of genus $g$. We
denote by $j$ a complex structure on $\Sigma$ and denote by $\CM_g
= \CM(\Sigma)$ the moduli space of complex structures on $\Sigma$.
We call a pair $((\Sigma,j),u)$ a $J$-holomorphic map if $u$ is
$(j,J)$-holomorphic, i.e., if it satisfies
$$
J\circ du = du \circ j.
$$
We say that $((\Sigma,j),u)$ is Fredholm regular if the
linearization of the map
$$
\delbar_J : (j,u) \mapsto \frac{du + J \circ du \circ j}{2}
$$
is surjective at $((\Sigma,j),u)$. We have the index formula  given
by \be\label{eq:index} \operatorname{Index} D_{(j,u)}\delbar_J =
\begin{cases} 2(c_1(M,\omega)(\beta) + (n-3)(1-g)) \quad & \mbox{for }\, g\geq 2 \\
2(c_1(M,\omega)(\beta) + 1) \quad & \mbox{for }\, g = 1 \\
2(c_1(M,\omega)(\beta) + n) \quad & \mbox{for }\, g = 0
\end{cases}
\ee for the maps $u$ with $[u] = \beta\in H_2(M,\Z)$, and hence
the virtual dimension of the associated moduli space
$\CM_g(M,J;\beta)$ is given by
\be\label{vdim}%
2(c_1(M,\omega)(\beta) + (n-3)(1-g))
\ee%
in all cases.

In this paper, we prove the following theorem.

\medskip

\begin{thm}\label{thm:embedded} Assume $n \geq 3$.
Let $(\Sigma,j)$ be a closed smooth Riemann surface of any genus
$g$. Then there exists a subset $\CJ_\omega^{emb} \subset
\CJ_\omega$ of second category such that for $J \in
\CJ_\omega^{emb}$, for any complex structure $j$ on $\Sigma$, all
somewhere injective $(j,J)$-holomorphic maps $u: \Sigma \to M$ are
Fredholm regular and \emph{embedded} whenever $c_1(M,\omega)([u])
\leq 0$.
\end{thm}

We would like to point out that for $n > 3$ this theorem has any
content only for the case $g = 0, \, 1$ for otherwise the
dimension formula \eqref{vdim} shows that higher genus somewhere
injective $J$-holomorphic curves cannot exist for a generic $J$
when $n > 3$.

In fact, we prove the following general theorem which
immediately gives rise to Theorem \ref{thm:embedded}.
We prove this theorem by establishing a transversality result for
the 1-jet evaluation map (see Proposition \ref{trans-Upsilon}) and then by
a dimension counting argument.

\begin{thm}\label{thm:main} Let $(M^{2n},\omega)$ be any
symplectic manifold and $\beta \in H_2(M,\Z)$. There exists a
subset $\CJ_\omega^{imm} \subset \CJ_\omega$ of second category
such that for $J \in \CJ_\omega^{imm}$, all somewhere injective
$(j,J)$-holomorphic maps $u:(\Sigma,j) \to (M,J)$ in class $\beta$
are immersed for any $j \in \CM_g$, provided \be c_1(\beta) +
(3-n)(g-1)< n -1. \ee And there exists another subset
$\CJ_\omega^{emb} \subset \CJ_\omega^{imm} \subset \CJ_\omega$ of
second category such that for $J \in \CJ_\omega^{emb}$ all
somewhere injective curves $(j,J)$-holomorphic curves are
embedded, provided \be c_1(\beta) + (3-n)(g-1)< n - 2. \ee
\end{thm}

An immediate corollary of Theorem \ref{thm:embedded}
is the following specialization to the
symplectic Calabi-Yau manifolds with $n \geq 3$.

\begin{cor} Let $(M,\omega)$ be symplectic Calabi-Yau, i.e.,
$(M,\omega)$ symplectic and  $c_1(M,\omega) = 0$. Assume $n \geq
3$. Let $(\Sigma,j)$ be a closed smooth Riemann surface of any
genus $g$. Then there exists a subset $\CJ_\omega^{emb} \subset
\CJ_\omega$ of second category such that for $J \in
\CJ_\omega^{emb}$, for any complex structure $j$ on $\Sigma$, all
somewhere injective $J$-holomorphic map $u: (\Sigma,j) \to (M,J)$
are Fredholm regular and \emph{embedded}.
\end{cor}

An immediate consequence of this corollary is the following
classification result of stable maps in
symplectic Calabi-Yau threefolds.

\begin{thm}\label{thm:stablemap} Suppose $c_1(M,\omega) = 0$ and $n = 3$.
Then there exists a subset $\CJ_\omega^{nodal} \subset \CJ_\omega^{emb} \subset \CJ_\omega$
of second category such that for $J \in \CJ_\omega^{nodal}$
any stable $J$-holomorphic map in Calabi-Yau
threefolds is one of the following three types :
\begin{enumerate}
\item it is either smooth and embedded or
\item it has smooth domain and factors through the composition
$$
u = u' \circ \phi : \Sigma \to C \hookrightarrow M
$$
for some embedding $u': C \to M$ and a ramified covering $\phi: \Sigma \to C$ or
\item it is a stable map of the type such that all of its irreducible
components have the same locus of images as $C$, an embedded
$J$-holomorphic curve in $M$, and are ramified over the domain of
$C$, except those of constant components.
\end{enumerate}
The latter two cases can occur only when $\beta = [u]$ is of the
form $\beta = d\gamma$ for some positive integer $d$ and homology
class $\gamma \in H_2(M,\Z)$.
\end{thm}

\medskip

A brief outline of this paper is in order.

In section 2, we establish the generic immersion property of
$J$-holomorphic curves (Theorem \ref{thm:immersion}) and prove the
first half of Theorem \ref{thm:main}. This section is the most novel
and essential part of this paper which involves the 1-jet evaluation
of the map $u$. Consideration of the 1-jet evaluation map in turn
forces us to work with the Fredholm setting of $W^{k,p}$ for $k \geq
3$ \emph{so that the 1-jet evaluation map becomes differentiable}
with respect to the variation of evaluation points, which involves
taking \emph{two} derivatives of the map. With this choice of
Sobolev spaces however, the actual Fredholm analysis involving the one-jet
evaluation map is rather delicate partly because one has to overcome
the fact that the evaluation of an $L^p$-map is not defined
pointwise. One novelty of our proof is a judicious usage of the
structure theorem of distributions with point supports. See the
proof of Lemma \ref{lem:eta=0}.

The scheme of our proof is motivated by a similar theorem of the
authors in \cite{oh-zhu} which sketched the proof of immersion
property of nodal Floer trajectory curves. The latter is in turn
partly motivated by Hutchings and Taubes' proof of Theorem 4.1
\cite{hutch-taubes} which concerns the immersion property of the
case of $4\,  (n=2)$ dimension in a different context : More
specifically see the proof of Lemma 4.2 \cite{hutch-taubes}. There
is also a Corollary 3.17 in \cite{wendl} which also concerns
immersion property of somewhere injective $J$-holomorphic curves for
the moduli spaces of dimension $0$ and $1$. In their proofs, both
papers utilize the fact that they concern a low dimensional moduli
space of $J$-holomorphic curves. One comparison between Theorem
\ref{thm:main} and Corollary 3.17 \cite{wendl} (with the case $\del
\Sigma = \emptyset$) is that the condition in our theorem is optimal
and corresponds to
$$
\operatorname{Index} D_{(j,u)}\delbar_J < 2(n-1)
$$
(for $n \geq 2$) while Wendl's would correspond to
$$
\operatorname{Index} D_{(j,u)}\delbar_J < 2.
$$
It is, however, conceivable that their proofs, with some
modifications, could be generalized to higher dimensional moduli spaces, which we did
not check.

Aside from establishing the immersion property, our \emph{natural}
Fredholm framework for the proof of 1-jet evaluation transversality
used in section \ref{sec:immersion} and \ref{appendix} has its own
merit and suits well for the generalization to the study of higher
jet evaluation transversality. We hope to come back to the study of
this higher order transversality elsewhere.

In section 3, we establish generic one-one property and prove
the second half of Theorem \ref{thm:main}.  Theorem
\ref{thm:embedded} stated above then follows by a dimension counting argument.

In section 4, we prove some key lemma, a type of removable
singularity theorem for which we utilize a structure theorem of distributions
with point supports (see \cite{gelfand} for example).

In section 5, we discuss some
implication of our results to the Gromov-Witten theory of Calabi-Yau
threefolds and derive Theorem \ref{thm:stablemap}.

The senior author would like to thank Pandharipande for his
interest on the current work and for some useful discussion on BPS
counts. We like to thank an anonymous referee for pointing out
some inaccuracies in our Fredholm setting of the previous version of
the paper, providing useful suggestions on improvement and pointing out
the reference \cite{wendl}.

\section{The 1-jet evaluation transversality}
\label{sec:immersion}

In this subsection, we will give the proof of immersion property.
This is the most novel and essential part of the present work. Except
the proof of this immersion property, the arguments used in other parts
are all standard and well-known in the study of pseudo-holomorphic curves.

\subsection{Fredholm setting}

We first provide some informal discussion to motivate the
necessary Fredholm set-up for the study of immersion property. We will
provide the precise analytical framework in the end of this discussion.

We consider a triple $(J,(j,u),z)$ of compatible $J$
and $u : (\Sigma,j) \to (M,J)$ a $(j,J)$-holomorphic map and $z \in \Sigma$.
Define a map $\Upsilon$ by
\be
\Upsilon(J,(j,u),z) = (\delbar(J,(j,u));\del(J,(j,u))(z))
\ee
where we denote
\beastar
\delbar(J,(j,u)) & : = & \delbar_{(j,J)}(u) = (du)^{(0,1)}_{(j,J)} = \frac{du + J du j}{2}\\
\del(J,(j,u)) &: = & \del_{(j,J)}(u) = (du)^{(1,0)}_{(j,J)} = \frac{du
- J du j}{2} . \eeastar
We now identify the domain and the target
of the map $\Upsilon$. For any given $(j,J)$, consider the bundles
over $\Sigma \times M$
\beastar
H^{(0,1)}_{(j,J)}(\Sigma \times M) &: = & \bigcup_{(z,x)} Hom''_{(j_z,J_x)}(T_z\Sigma,T_x M) \\
H^{(1,0)}_{(j,J)}(\Sigma \times M) &: = & \bigcup_{(z,x)}
Hom'_{(j_z,J_x)}(T_z\Sigma,T_x M), \eeastar where the above unions
are taken for all $(z,x)$ of $(\S\times M)$. Over any $(z,x)$, the
fibers are the $(j,J)$-anti-linear and $(j,J)$-linear parts of
$Hom(T_z\Sigma,T_xM)$, denoted by $Hom^{''}_{(j_z,J_x)}(T_z \S,T_x
M)$ and $Hom^{'}_{(j_z,J_x)}(T_z \S,T_x M)$ respectively. We denote
$$
\Lambda_{(j_z,J_x)}^{(1,0)}(T_x M) = Hom'_{(j_z,J_x)}(T_z\Sigma,T_x M)
$$
as usual.

We now introduce the necessary framework for the
Fredholm theory needed to prove the main theorem.
Let $\beta \in H_2(M,\Z)$ be given and consider
the off-shell function space
$$
\CF(\Sigma,M;\beta) = \{((\Sigma,j),u) \mid j \in \CM(\Sigma),\,
u: \Sigma \to M, [u]=\beta\}
$$
hosting the operator $\delbar_J:(j,u) \mapsto \delbar_{(j,J)}(u)$.

For each given $(J,(j,u),z)$, we associate a $2n$-dimensional vector
space
$$
H^{(1,0)}_{(J,(j,u),z)} : = \Lambda_{(j,J)}^{(1,0)}(u^*TM)|_{z} =
\Lambda_{(j_z,J_{u(z)})}^{(1,0)}(T_{u(z)}M)
$$
and define the vector bundle of rank $2n$
$$
H^{(1,0)} = \bigcup_{(J,(j,u),z)}\Lambda_{(j,J)}^{(1,0)}(u^*TM)|_{z}
$$
over the space $\CF_1(\Sigma,M;\beta)$ defined by
$$
\CF_1(\Sigma,M;\beta) = \{((\Sigma,j),u,z) \mid ((\Sigma,j),u) \in
\CF(\Sigma,M;\beta),\, z \in \Sigma\}.
$$
We denote the corresponding moduli space of marked $J$-holomorphic
curves $((\Sigma,j),u,z)$ by $\CM_1(\Sigma,M;\beta)$.

\begin{rem}
In this paper, the domain complex structure $j$ does not play much
role in our study. Especially it does not play any role throughout
our calculations except that it appears as a parameter.
\end{rem}

We introduce the standard bundle
$$
\CH'' = \bigcup_{(J,(j,u))} \CH''_{((j,u),J)}, \quad \CH''_{(J,(j,u))} = \Omega_{(j,J)}^{(0,1)}(u^*TM).
$$
Then we have the map
$$
\Upsilon : \CJ_\omega \times \CF_1(\Sigma,M;\beta) \to \CH''\times
H^{(1,0)} ;\, (J,(j,u),z) \mapsto (\delbar_{(j,J)}
u,\;(\del_{(j,J)}u)(z))
$$
where $\CH''\times H^{(1,0)}$ is the fiber product
of the two bundles
$$
\pi_1: \CH'' \to \CJ_\omega \times \CF(\Sigma,M;\beta)
$$
and
$$
\pi_2 : H^{(1,0)} \to \CJ_\omega  \times \CF_1(\Sigma,M;\beta)
\to \CJ_\omega \times \CF(\Sigma,M;\beta).
$$
More explicitly we can express the fiber product as
$$
\CH''\times H^{(1,0)} := \left\{(\eta, \zeta_0; J,(j,u),z) \, \Big| \,
\eta \in \CH''_{(J,(j,u))}, \, \zeta_0 \in H^{(1,0)}_{(J,(j,u),z)} \right\}
$$
We regard this fiber product as a vector bundle over $\CJ_\omega \times \CF_1(\Sigma,M;\beta)$,
$$
(\eta, \zeta_0; J,(j,u),z) \mapsto (J,(j,u),z)
$$
whose fiber at $(J,(j,u),z)$ is given by
$$
\CH''_{(J,(j,u))} \times  H^{(1,0)}_{(J,(j,u),z)}.
$$
Then the above map $\Upsilon$ will become a smooth
\emph{section} of this vector bundle.

The union of standard moduli spaces $\CM_1(M,J;\beta)$ over $J \in \CJ_\omega$
is nothing but
\be
\Upsilon^{-1}(o_{\CH''} \times H^{(1,0)})/\operatorname{Aut}(\Sigma)
\ee
where $o_{\CH''}$ is the
zero section of the bundle $\CH''$ defined above, and
$\operatorname{Aut}(\Sigma)$ acts on $((\S,j),u)$ by conformal
equivalence for any $j$.  We also denote
 \beastar
\widetilde \CM_1(M;\beta) & = & \Upsilon^{-1}(o_{\CH''} \times H^{(1,0)})\\
\widetilde \CM_1(M,J;\beta) & = & \widetilde \CM_1(M;\beta) \cap
\pi_2^{-1}(J). \eeastar

\medskip

The following characterization of the critical point is obvious to
see, which however is a key ingredient for the Fredholm framework
used in our proof of immersion property.

\begin{lem}
For any $((j,u),z) \in \widetilde \CM_1(M,J;\beta)$, since
$\delbar_{J,j}u=0$, we have \be \label{eq:key} du(z) = 0 \quad
\mbox{if and only if} \quad \del_{(j,J)}u(z) = 0.
\ee
\end{lem}

Some remarks concerning the necessary Banach manifold set-up of the map
$\Upsilon$ are now in order :
\begin{enumerate}
\item To make evaluating $\del u$ at a point $z \in
\Sigma$ make sense, we need to take at least $W^{2,p}$-completion
with $p > 2$ of $\CF(\S,M;\b)$ so that $\delbar_{(j,J)} u$ lies in
$W^{1,p}$ which is then continuous. We actually need to take
$W^{k,p}$-completion of $\CF(\S,M;\b)$ with $k \geq 3$ so that the section
$\Upsilon$, especially the evaluation map, is differentiable (see \eqref{eq:DUpsilon}).
We denote the corresponding completion of $\CF(\Sigma,M;\beta)$ by
$$
\CF^{k,p} = \CF^{k,p}(\Sigma,M;\beta).
$$
\item We provide the $\CH''$ the topology of a $W^{k-1,p}$ Banach bundle,
with each fiber of class $W^{k-1,p}$. The choice of $k$ will also
depend on the index of the linearization of $D\delbar_{(j,J)}$
on $\CF(\Sigma,M;\beta)$ and should be chosen sufficiently large so
that one can apply Sard-Smale theorem \cite{smale}.
\item We also need to provide some Banach manifold structure on $\CJ_\omega$.
We can borrow Floer's scheme \cite{floer} for this whose details we
refer readers thereto. Also see Remark 3.2.7 \cite{mcd-sal04}.
\end{enumerate}

We will assume these settings during the proof of Proposition
\ref{trans-Upsilon} without explicit mentioning unless it is
absolutely necessary. At fixed $(J, (u,j),z_0)$ where we do linearization of $\Upsilon$,
we will write
\beastar
\Omega^0_{k,p}(u^*TM) & := & W^{k,p}(u^*TM) = T_u \CF^{k,p}(\Sigma,M;\beta)\\
\Omega^{(0,1)}_{k-1,p}(u^*TM) & := & W^{k-1,p}\left(\Lambda_{(j,J)}^{(0,1)}(u^*TM)\right)
\eeastar
for the simplicity of notations. Let $o_{H^{(1,0)}}$ be
the zero section of $H^{(1,0)}$.

\subsection{Proof of generic immersion property}

We now prove the following proposition by linearizing the section
$\Upsilon$.

\begin{prop}\label{trans-Upsilon} The section $\Upsilon$ is transverse to the
zero section
\be o_{\CH''\times H^{(1,0)}} = o_{\CH''} \times o_{H^{(1,0)}} \subset \CH'' \times
H^{(1,0)}.
\ee
In particular the set
$$
\Upsilon^{-1}(o_{\CH''} \times o_{H^{(1,0)}})
$$
is a submanifold of $\widetilde \CM_1(M;\beta)$ of codimension $2n$.
\end{prop}

\begin{proof} Recall that the subset
\be\label{eq:oCH''}%
o_{\CH''} \times o_{H^{(1,0)}}\subset o_{\CH''}\times
H^{(1,0)}
\ee%
is a submanifold of codimension $2n$. So it is easy to check the
statement on the codimension once we prove $\Upsilon$ is transverse
to the submanifold $ o_{\CH''} \times o_{H^{(1,0)}} \subset \CH''
\times H^{(1,0)}$.

Let $(J,(j,u),z) \in \Upsilon^{-1}(o_{\CH''} \times o_{H^{(1,0}})$.
Pick any $J$-complex connection, i.e., $\nabla$ with $\nabla J = 0$
and denote by $\nabla_{du}$ the pull-back connection of $\nabla$ by
$u$.

The linearization of $\Upsilon$ at $(J,(j,u),z)$
$$
D_{(J,(j,u),z)}\Upsilon : T_J\CJ_\omega \times T_{((j,u),z)}\CF_1^{k,p}(\Sigma,M;\beta) \to
\CH''_{(J,(j,u))} \times H^{(1,0)}_{(J,(j,u),z)}
$$
is given by the formula
\be\label{eq:DUpsilon}%
(B,(b,\xi),v)) \mapsto
\Big(D_{J,(j,u)}\delbar(B,(b,\xi)),\;D_{J,(j,u)}\del(B,(b,\xi))(z) +
\nabla_{du(v)} (\del_{(j,J)}u) \Big)
\ee%
for $B\in T_J\CJ_{\o}, b\in T_j\CM(\S), v\in T_z\S$ and $\xi \in
\Omega^0_{k,p}(u^*TM)$. Recall that $u$ is in $W^{k,p}$ with $k \geq
3$ (in fact, $u$ is smooth by elliptic regularity since
$\delbar_{(j,J)}u = 0$) so $D_{J,(j,u)}\del(B,(b,\xi))$ and
$\nabla_{du(v)} (\del_{(j,J)}u)$ are in $W^{k-2,p}$ where $k-2
\geq1$. Therefore their evaluations at $z$ are well-defined.

We need to prove that at each $(J,(j,u),z_0) \in
\Upsilon^{-1}(o_{\CH''} \times o_{H^{(1,0)}})$, the system of
equations
\bea%
D_{J,(j,u)}\delbar(B,(b,\xi)) &=& \gamma \label{eq:Dbar}\\
D_{J,(j,u)}\del(B,(b,\xi))(z_0) + \nabla_{du(v)} (\del_{(j,J)}u)) &
= & \zeta_0 \label{eq:DJueta-}
\eea%
has a solution $(B,(b,\xi),v)$ for each given data
$$
\gamma \in \Omega^{(0,1)}_{k-1,p}(u^*TM), \quad \zeta_0 \in H^{(1,0)}_{(J,(j,u),z_0)}.
$$
It will be enough to consider the triple with $b = 0$ and $v=0$ which we will assume
from now on.

In general, a well-known computation shows \be\label{eq:formula}
D_{J,(j,u)}\del(B,(0,\xi)) = (\nabla_{du} \xi)^{(1,0)}_{(j,J)} +
T^{(1,0)}_{(j,J)}(du, \xi)+\frac{1}{2}B \circ du \circ j \ee with
respect to a $J$-complex connection $\nabla$ and its torsion tensor
$T$. Here we denote
$$
T^{(1,0)}_{(j,J)}(du,\xi) = \frac{1}{2}\left(T(du,\xi) + J T(du \circ j, \xi)\right).
$$
However if $ u\in
\Upsilon^{-1}(o_{\CH''} \times o_{H^{(1,0)}_{(j,J)}})$, we have $du(z_0)
= 0$ and hence
$$
T^{(1,0)}_{(j,J)}(du(z_0), \xi(z_0))=0 =\frac{1}{2}B(u(z_0))
\circ du(z_0) \circ j_{z_0}
$$
for any $\xi$. If we just want to solve
\eqref{eq:DJueta-} at $z_0$, then \eqref{eq:DJueta-} is reduced to
\be \label{eq:Dpt}%
(\nabla_{du} \xi)^{(1,0)}_{(j,J)}(z_0)=\zeta_0.
\ee%

Now we study solvability of \eqref{eq:Dbar}-\eqref{eq:DJueta-} by applying the
Fredholm alternative. For this purpose, we make the following crucial remark

\begin{rem} \label{rem:3pversus2p}
We emphasize  that for the map
$$
(z_0,v) \mapsto D_{J,(j,u)}\del(B,(b,\xi))(z_0) + \nabla_{du(v)}
(\del_{(j,J)}u))
$$
to be defined as a continuous map to $H^{(1,0)}_{(J,(j,u),z_0)} = \Lambda^{(1,0)}(T_{u(z_0)}M)$, the map $u$ must
be at least $W^{2 +\epsilon, p}$ for $\e > 0$ : On $W^{2,p}$, the map
$D_{J,(j,u)}\del(B,(b,\xi))$ will be only in $L^p$ for which
the evaluation at a point is not defined in general, let alone being continuous.
However, the evaluation map
\be\label{eq:W2ptoTM}
z_0 \mapsto D_{J,(j,u)}\del(B,(0,\xi))(z_0)
\ee
is well-defined and continuous on $W^{2,p}$ as shown by the explicit formula
\eqref{eq:formula}, which involves only \emph{one} derivative of the section $\xi$.
This reduction from $W^{k,p}$ to $W^{2,p}$ of the regularity requirement
in the study of the map \eqref{eq:W2ptoTM}, which can be achieved after restricting to $b=0$, $v = 0$,
will play a crucial role in our proof. See the proof of Lemma \ref{lem:eta=0}.
\end{rem}

Utilizing this remark, we will first show that the image of the map
\eqref{eq:DUpsilon} restricted to the elements of the form
$(B,(0,\xi),0)$ is onto as a map
$$
T_J \JJ_\omega \times \Omega^0_{2,p}(u^*TM) \to \Omega^{(0,1)}_{1,p}(u^*TM) \times
H^{(1,0)}_{(J,(j,u),z_0)}
$$
where $(u,j,z_0,J)$ lies in $o_{\CH^{''}}\times o_{H^{(1,0)}}$. In the end of
the proof, we will establish solvability of \eqref{eq:Dbar}-\eqref{eq:DJueta-}
on $W^{k,p}$ for $\gamma \in W^{k-1,p}$ by applying an elliptic regularity result
of the map \eqref{eq:Dbar}.

We regard
$$
\Omega^{(0,1)}_{1,p}(u^*TM) \times H^{(1,0)}_{(J,(j,u),z_0)}: = \CB
$$
as a Banach space with the norm
$$
\|\cdot\|_{1,p} + |\cdot|
$$
where $|\cdot|$ any norm induced by an inner product on
$$
H^{(1,0
)}_{(J,(j,u),z_0)}=\Lambda^{(1,0)}_{(j,J)}(u^*TM)_{z_0} \cong \C^n.
$$
For the clarification of notations, we denote the natural pairing
$$
\Omega^{(0,1)}_{1,p}(u^*TM)\times \left(\Omega^{(0,1)}_{1,p}(u^*TM)\right)^*\to \R
$$
by $\langle \cdot, \cdot
\rangle$ and the inner product on $H^{(1,0)}_{(J,(j,u),z_0)}$ by $(\cdot,
\cdot)_{z_0}$.

We will first prove that the image is dense in $\CB$.
\par
Let $(\eta, \alpha_{z_0})
\in\left(\Omega^{(0,1)}_{1,p}(u^*TM)\right)^* \times
H^{(1,0)}_{(J,(j,u),z_0)}$ such that \be\label{eq:0} \langle
D_{J,(j,u)}\delbar_{(j,J)}(B,(0,\xi)), \eta \rangle +
(D_{J,(j,u)}\del_{(j,J)}(B,(0,\xi))(z_0), \alpha_{z_0})_{z_0}= 0
\ee for all $\xi \in \Omega^0_{2,p}(u^*TM)$ and $B\in T_J\CJ_\o$.
Without loss of any generality, we may assume that $\xi$ is smooth
since $C^\infty(u^*TM) \hookrightarrow \Omega^0_{2,p}(u^*TM)$ is
dense. Under this assumption, we would like to show that $\eta = 0
= \alpha_{z_0}$.

Now we simplify the expression of $D_{J,(j,u)}\del_{(j,J)}(B,(0,\xi))(z_0)$
in complex coordinates $z$ at $z_0$.
Let $x_0=u(z_0)$, and identify a
neighborhood of $z_0$ with an open subset of $\C$ and a neighborhood
of $x_0$ with an open set in $T_{x_0}M$. We now
introduce the linear operator $q_{J,x_0}$ defined by
$$
q_{J,x_0}(x) = (J_{x_0} + J(x))^{-1}(J_{x_0} - J(x))
$$
for $x$ such that $d(x,x_0) < \delta$ for $\delta > 0$ depending
only on $(M,\omega,J)$ but independent of $x_0$. $q_{J,x_0}$
satisfies $q_{J,x_0}(x_0) = 0$. (See \cite{sikorav}.) Then if we
identify $(T_{x_0} M,J_{x_0}) \cong \C^n$, we can write the
operator
\beastar%
(\nabla_{du} \xi)^{(1,0)}_{(j,J)} &=& \del \xi - q_{J,x_0}(u)\delbar \xi +
C\cdot \xi\\
&=& \del \xi - A\cdot\delbar \xi + C\cdot \xi,
\eeastar%
where in a neighborhood of $z_0$, $\del, \, \delbar$ are the
standard Cauchy-Riemann operators on $\C^n$ and $A, \, C$ are smooth pointwise
(matrix) multiplication operators with \be\label{eq:ABr-} A(z_0) =
C(z_0) = 0.\ee
Therefore we have
\be\label{DJju}%
D_{J,(j,u)}\del_{(j,J)}(B,(0,\xi))(z_0) = (\nabla_{du}
\xi)^{(1,0)}_{(j,J)}(z_0) = (\del \xi - A \cdot \delbar \xi + C
\cdot \xi)(z_0)
\ee%
at the given point $z_0$ for any given $\zeta_0$. Since
we just need $\xi$ to satisfy \eqref{DJju} at $z_0$, by the condition
of $A$ and $C$ at $z_0$, we have shown
\be\label{pure-d}%
D_{J,(j,u)}\del_{(j,J)}(B,(0,\xi))(z_0) = \del \xi(z_0).
\ee%

By the above discussion on $D_{J,(j,u)}\delbar_{(j,J)}(B,(0,\xi))$
and $D_{J,(j,u)}\del_{(j,J)}(B,(0,\xi))(z_0)$, \eqref{eq:0} is equivalent to
\be\label{eq:Tgt}%
\langle D_u\delbar_{(j,J)} \xi+ \frac{1}{2}B\circ du\circ j, \eta \rangle
+ \langle \del \xi, \delta_{z_0} \alpha_{z_0}\rangle =0
\ee%
for all $B$ and $\xi$ of $C^{\infty}$ where $\delta_{z_0}$ is the
Dirac-delta function.

Taking $B=0$ in \eqref{eq:Tgt}, we obtain
\be \label{eq:coker} %
\langle D_u\delbar_{(j,J)} \xi, \eta \rangle + \langle \del \xi,
\delta_{z_0} \alpha_{z_0} \rangle = 0 \quad \mbox{ for all $\xi$
of $C^{\infty}$ }.
\ee %
Therefore by definition of the distribution derivatives, $\eta$ satisfies
$$
(D_u\delbar_{(j,J)})^\dagger \eta - \delbar(\delta_{z_0} \alpha_{z_0}) = 0
$$
as a distribution, i.e.,
$$
(D_u\delbar_{(j,J)})^\dagger \eta = \delbar(\delta_{z_0} \alpha_{z_0})
$$
where $(D_u \delbar_{(j,J)})^\dagger$ is the formal adjoint of $D_u
\delbar_{(j,J)}$ whose symbol is the same as $D_u\del_{(j,J)}$ and so is
an elliptic first order differential operator. We also recall that
$\del^\dagger = - \delbar$. Since $\supp \delbar(\delta_{z_0})
\alpha_{z_0} \subset \{z_0\}$, we have $(D_u\delbar_{(j,J)})^\dagger
\eta = 0$ on $\Sigma \setminus \{z_0\}$ as a distribution. Then by
the elliptic regularity (see Theorem 13.4.1 \cite{Ho}, for example),
$\eta$ must be smooth on $\Sigma \setminus \{z_0\}$.

On the other hand, by setting $\xi = 0$ in \eqref{eq:Tgt},  we get
\be
\label{eq:1/2B} \langle \frac{1}{2}B\circ du\circ j, \eta
\rangle=0
\ee
for all $B\in T_{J} \CJ_\omega$. From this identity,
standard argument from \cite{floer}, \cite{mcduff} shows that $\eta=0$ in a
small neighborhood of any somewhere injective point in $\Sigma
\setminus \{z_0\}$. Such a somewhere injective point exists by the
hypothesis of $u$ being somewhere injective and the fact that the
set of somewhere injective points is open and dense in the domain
under the hypothesis (see \cite{mcduff}). Then by the unique
continuation theorem, we conclude that $\eta = 0$ on $ \Sigma
\backslash\{z_0\}$ and so the support of $\eta$ as a distribution
on $\Sigma$ is contained at the one-point subset $\{z_0\}$ of
$\Sigma$.

We will postpone the proof of the following lemma till section
\ref{appendix}, in which reduction of the regularity requirement
for $u$ mentioned in Remark \ref{rem:3pversus2p} will play an essential role.

\begin{lem} \label{lem:eta=0}  $\eta$ is a distributional solution of
$(D_u \delbar_{(j,J)})^\dagger \eta = 0$ on $\Sigma$
and so continuous.
In particular, we have $\eta = 0$ in $\left(\Omega^{(0,1)}_{(1,p)}(u^*TM)\right)^*$.
\end{lem}

Once we know $\eta = 0$, the equation \eqref{eq:0} is reduced to
\be\label{eq:simple0}
(D_{J,(j,u)}\del_{(j,J)}(B,(0,\xi))(z_0), \alpha_{z_0})_{z_0} = 0
\ee
It remains to show that $\alpha_{z_0} = 0$. For this, we have only to
show that the image of the evaluation map
$$
\xi \mapsto D_{J,(j,u)}\del_{(j,J)}(0,(0,\xi))(z_0) = \del \xi(z_0)
$$
is surjective onto
$H^{(1,0)}_{(J,(j,u),z_0)}=\Lambda_{(j,J)}^{(1,0)}(T_{u(z_0)}M)$. The equality
comes from \eqref{pure-d}.

To show this surjectivity, we need to prove the existence of $\xi$ satisfying
\be\label{ac}%
\del \xi (z_0) = \zeta_0
\ee%
at the given point $z_0$ for any given $\zeta_0$.
We can multiply a cut-off function $\chi$ to $\zeta_0$ with
$\chi \equiv 1$ to make $\zeta(z): = \chi(z) \zeta_0$ supported in a
sufficient small neighborhood around $z_0$, and apply Cauchy
integral formula in coordinates to solve
$$
\del \xi = \zeta
$$
in some neighborhood around $z_0$. This finishes the existence of a
solution to (\ref{eq:Dpt}) and hence \eqref{eq:DJueta-} and so
the proof of the claim that the image of \eqref{eq:DUpsilon}
with $v=0$ is \emph{dense} in
$$
\Omega^{(0,1)}_{1,p}(u^*TM) \times H^{(1,0)}_{(J,(j,u),z_0)}.
$$
We recall that for any fixed $(J,j)$, the
image of $D_u\delbar_{(j,J)}$ on $T_J\JJ_\omega \times \Omega^0_{2,p}(u^*TM)$ is
closed in $\Omega^{(0,1)}_{(1,p)}(u^*TM)$ and that $H^{(1,0)}_{(J,(j,u),z_0)}$ is a
finite dimensional vector space. Therefore the image of the linearization
$\eqref{eq:DUpsilon}$ is also \emph{closed}. Hence $\eqref{eq:DUpsilon}$
is surjective
$\Omega^{(0,1)}_{(1,p)}(u^*TM) \times H^{(1,0)}_{(J,(j,u),z_0)}$
as a map from $T_J\JJ_\omega \times T_{((j,u),z)}\CF_1^{2,p}(\Sigma,M;\beta)$.

Now finally suppose $\gamma$ is in the subspace $\Omega^{(0,1)}_{k-1,p}(u^*TM) \subset
\Omega^{(0,1)}_{1,p}(u^*TM)$ of higher regularity.
We recall $k \geq 3$. By the above analysis of the map
$\eqref{eq:DUpsilon}$ for $b = 0 = v$, we can find a solution
$(B,(0,\xi),0)$ of \eqref{eq:Dbar}-\eqref{eq:DJueta-} with $B \in
T_J\JJ_\omega$ and with \emph{$\xi$ as an element in $\Omega^0_{2,p}(u^*TM) = W^{2,p}(u^*TM)$}
for any $\zeta_0 \in H^{(1,0)}_{(J,(j,u),z_0)}$. By elliptic regularity of
\eqref{eq:Dbar}, $\xi$ indeed lies in $W^{k,p}$ if $\gamma
\in W^{k-1,p}$, and hence in $\Omega^0_{k,p}(u^*TM)$.
Therefore the map $\eqref{eq:DUpsilon}$ is onto, i.e., $\Upsilon$ is
transverse to the submanifold $ o_{\CH''} \times
o_{H^{(1,0)}} \subset \CH'' \times H^{(1,0)}$.
\end{proof}

Finally we have the natural projection
$$
\pi:\widetilde \CM_1(M;\beta) : = \bigcup_{J \in \CJ_\omega}
\widetilde \CM_1(M,J;\beta) \to \CJ_\omega.
$$
The projection has index $2(c_1(\beta) + n(1-g)) + 2$, so for any
regular value $J$, the moduli space
$$
\widetilde\CM_1^{crit}(M,J;\b):=\Upsilon^{-1} (o_{\CH''} \times
o_{H^{(1,0)}} )\cap \pi^{-1}(J)
$$
is of dimension
$$
2(c_1(\beta) + n(1-g)) + 2-2n.
$$
Let the moduli space
$$
\CM_1^{crit}(M,J;\b):=\widetilde\CM_1^{crit}(M,J;\b)/\Aut(\Sigma),
$$
where $\Aut(\Sigma)$ acts on marked Riemann surfaces $((\S,j),z)$ by
conformal equivalence then on the maps from them. Geometrically
$\CM_1^{crit}(M,J;\b)$ consists of $J$-holomorphic curves in class
$\b$ with at least one critical point. As a smooth orbifold, we have
$$
\text{dim }\CM_1^{crit}(M,J;\b)=2 \big(c_1(\b)+(3-n)(g-1)+1-n\big).
$$
Therefore, $\CM_1^{crit}(M,J;\b)$ is empty whenever this dimension is negative, i.e.,
$$
c_1(\b)+(3-n)(g-1)<n-1.
$$

We just set
$$
\CJ_\omega^{imm} = \mbox{the set of regular values of $\pi$}
$$
which finishes the proof of the following theorem.

\begin{thm} \label{thm:immersion} There exists a subset $\CJ_\omega^{imm} \subset
\CJ_\omega$ of second category such that for $J \in
\CJ_\omega^{imm}$ all somewhere injective $(j,J)$-holomorphic maps
$u:\Sigma \to M$ are immersed for any $j \in \CM_g$, provided
\be\label{eq:immersed} c_1(\beta) + (3-n)(g-1)< n -1 \ee
\end{thm}

\section{Proof of Lemma \ref{lem:eta=0}}
\label{appendix}

In this section, we prove Lemma \ref{lem:eta=0}. Our primary goal
is to prove \be\label{eq:tildexi} \langle D_u\delbar_{(j,J)} \xi,
\eta \rangle = 0 \ee for all smooth $\xi \in \Omega^0(u^*TM)$,
i.e., $\eta$ is a distributional solution of
$(D_u\delbar_{(j,J)})^\dagger \eta = 0$ \emph{on the whole
$\Sigma$}, not just on $\Sigma \setminus \{z_0\}$ which was shown
in section \ref{sec:immersion}.

We start with \eqref{eq:coker} \be\label{eq:coker-append} \langle
D_u\delbar_{(j,J)} \xi, \eta \rangle + \langle \del \xi,
\delta_{z_0} \alpha_{z_0} \rangle = 0 \quad \mbox{for all $\xi$ of
$C^{\infty}$ }. \ee We first simplify the expression of the
pairing $\langle D_u\delbar_{(j,J)} \xi, \eta \rangle$ knowing that
$\supp \eta \subset \{z_0\}$.

Let $z$ be a complex coordinate centered at a fixed marked point $z_0$
and  $(w_1, \cdots, w_n)$ be the complex coordinates on $T_{u(z_0)}M$
regarded as coordinates on a neighborhood of $u(z_0)$.
We consider the standard metric
$$
h = \frac{\sqrt{-1}}{2} dz d\bar z
$$
on a neighborhood $U$ of $z_0$ and with respect to the coordinates
$(w_1,...,w_n)$ we fix any Hermitian metric on $\C^n$.

The following lemma will be crucial in our proof.

\begin{lem}
For any smooth section $\xi$ of $u^*(TM)$ and $\eta$ of
$\left(\Omega^{(0,1)}_{1,p}(u^*TM)\right)^*$
$$
\langle D_u\delbar_{(j,J)}\xi , \eta \rangle = \langle \delbar \xi,
\eta \rangle,
$$
 where $\delbar$ is the standard Cauchy-Riemann operators on $\C^n$ in the above coordinate.
\end{lem}
\begin{proof} We have already shown that
$\eta$ is a distribution with $\supp \eta \subset \{z_0\}$. By the
structure theorem on the distribution supported at a point $z_0$
(see section 4.5, especially p. 119, of \cite{gelfand}, for
example), we have
$$
\eta = P\left(\frac{\del}{\del s}, \frac{\del}{\del t}\right)(\delta_{z_0})
$$
where $z = s + it$ is the given complex coordinates at $z_0$ and
$P\left(\frac{\del}{\del s}, \frac{\del}{\del t}\right)$ is a differential
operator associated by the polynomial $P$ of two variables with coefficients
in $\left(\Lambda^{(0,1)}_{(j_{z_0},J_{u(z_0)})}(u^*TM)\right)^*$.

Furthermore since $\eta \in (W^{1,p})^*$, the degree of $P$ \emph{must be zero}
and so we obtain
\be\label{eq:eta=adelta}
\eta = \alpha_{z_0} \delta_{z_0}
\ee
for some constant vector $\alpha_{z_0}$ : This is because the `evaluation
at a point of the derivative' of $W^{1,p}$ map does not
define a continuous functional on $W^{1,p}$.

We can write
$$
D_u\delbar_{(j,J)}\xi = \delbar \xi + E \cdot \del \xi + F\cdot \xi
$$
near $z_0$ in coordinates similarly as we did in \eqref{DJju} for
the operator $D_u\del_{(j,J)}$, where $E$ and $F$ are zero-order
matrix operators with $E(z_0) = 0 = F(z_0)$. Combining this with
\eqref{eq:eta=adelta}, we derive
$$
\langle E \cdot \del \xi + F \cdot \xi, \eta \rangle = \langle E \cdot \del \xi
+ F\cdot \xi, \alpha_{z_0} \delta_{z_0} \rangle
= (E(z_0) \del \xi(z_0) + F(z_0) \xi(z_0), \alpha_{z_0})_{z_0} = 0.
$$
Therefore we obtain
$$
\langle D_u\delbar_{(j,J)}\xi , \eta \rangle = \langle \delbar \xi + E \cdot \del \xi + F\cdot \xi,
\eta \rangle = \langle \delbar \xi, \eta \rangle
$$
which finishes the proof.
\end{proof}

\begin{rem}
We note that \eqref{eq:eta=adelta} is where we needed to have made the
reduction of the regularity requirement for the map
$\xi$ from $W^{k,p}$ to $W^{2,p}$ in the first half of the proof of
Proposition \ref{trans-Upsilon} : if we had not made the reduction but
required $\xi$ to be in $W^{k,p}$, its derivative would have lied in $W^{k-1,p}$ and
hence we could have only concluded
$$
\eta = P\left(\frac{\del}{\del s}, \frac{\del}{\del t}\right)(\delta_{z_0})
$$
with $\deg P \leq k-1$. This would then cause a problem in the argument
for the rest of the proof of this lemma.
\end{rem}

This lemma then implies that \eqref{eq:coker-append} becomes
\be\label{eq:coker-simple}
\langle \delbar \xi, \eta \rangle + \langle \del \xi, \delta_{z_0} \alpha_{z_0}
\rangle = 0 \quad
\mbox{for all $\xi$}.
\ee
We decompose $\xi$ as
$$
\xi(z) = \widetilde \xi(z) + \chi(z) (z-z_0) a(z_0)
$$
by defining $\widetilde\xi$ by
$$
\widetilde \xi(z) = \xi(z) - \chi(z) (z-z_0)a(z_0)
$$
where we express the one-form $\del \xi$ as \be\label{eq:delxiz}
\del \xi(z) = a(z) dz \ee on $U$ in coordinates with $a(z_0) \in
\C^n$, and $\chi$ is a cut-off function with $\chi \equiv 1$ in a
small neighborhood $V$ of $z_0$ and satisfies $\supp \chi \subset
U$. The choice of this decomposition is dictated by the fact
\be\label{eq:delchi} \del\left(\chi(z) (z-z_0) a(z_0)\right) (z_0)
= a(z_0) dz. \ee Then $\widetilde \xi$ is a smooth section on
$\Sigma$, and satisfies
$$
\quad \del \widetilde \xi(z_0) = 0, \quad
\delbar \widetilde \xi = \delbar{\xi} \quad \mbox{on $V$.}
$$
Therefore applying \eqref{eq:coker-simple} to $\widetilde \xi$ instead of $\xi$, we have
$$
\langle \delbar \widetilde \xi, \eta \rangle + \langle \del \widetilde \xi, \delta_{z_0} \alpha_{z_0}
\rangle = 0.
$$
But we have
\be\label{eq:tildexi=xi}
\langle \delbar \widetilde \xi, \eta \rangle = \langle \delbar \xi, \eta \rangle
\ee
since $\delbar \widetilde \xi = \delbar{\xi}$ on $V$ and $\supp \eta \subset \{z_0\}$.
Again using the support property $\supp \eta \subset \{z_0\}$
and \eqref{eq:delxiz}, \eqref{eq:delchi},  we derive
\bea\label{eq:delxi-adz}
\langle \del \widetilde \xi, \delta_{z_0} \alpha_{z_0} \rangle & = &
\langle \del \xi, \delta_{z_0} \alpha_{z_0} \rangle - \langle \del(\chi(z)(z-z_0)a(z_0)),
\delta_{z_0}\alpha_{z_0} \rangle \nonumber\\
& = & (\del \xi(z_0),\alpha_{z_0})_{z_0} - (a(z_0)dz, \alpha_{z_0})_{z_0} \nonumber\\
& = & (\del \xi(z_0)- a(z_0)dz,\alpha_{z_0})_{z_0} = 0.
\eea
Substituting \eqref{eq:tildexi=xi} and
\eqref{eq:delxi-adz} into \eqref{eq:coker-simple}, we obtain
$$
\langle \delbar \xi, \eta \rangle = 0
$$
and so we have finished the proof of \eqref{eq:tildexi}.

By the elliptic regularity, $\eta$ must be smooth. Since we have
already shown $\eta = 0$ on $\Sigma \setminus \{z_0\}$, continuity
of $\eta$ proves $\eta = 0$ on the whole $\Sigma$.

\section{Generic one-one property}
\label{sec:embedding}

In this subsection, we will give the embedding part of the main theorem.
As before we assume that $(j,u)$ is a somewhere injective $J$-holomorphic
curve in class $\beta \in H_2(M,\Z)$.

We note that for any given self-intersection point $x \in \mbox{Im } u$
we have
$$
u^{-1}(x) = \{z_1, \cdots, z_k\}
$$
for some integer $k \geq 2$. So we will consider the moduli
space of $J$-holomorphic maps with 2 marked
points and in the homology class $\beta$.

We consider the triples $(J,(j,u),(z_1, z_2))$ and the map
\be \Upsilon_2 :
(J,(j,u),(z_1,z_2)) \mapsto (\delbar_{(j,J)}u, (u(z_1),u(z_2))).
\ee
We introduce the necessary framework for the Fredholm theory
needed to prove Theorem \ref{thm:immersion}. Similarly as
$\CF_1(\Sigma,M;\beta)$, we define
$$
\CF_2(\Sigma,M;\beta) = \{ ((j,u),(z_1, z_2)) \mid (j,u) \in
\CF(\Sigma,M;\beta), \, z_1, \, z_2 \in \Sigma, z_1\neq z_2 \}
$$
and by $\CM_2(M,J;\beta) = \widetilde \CM_2(M,J;\beta)/\Aut(\Sigma)$ the corresponding moduli spaces of
$J$-holomorphic curves where
$$
\widetilde \CM_2(M,J;\beta) = \{ ((j,u),(z_1, z_2)) \in \CF_2(\Sigma,M;\beta) \mid
\delbar_{j,J}u = 0 \}.
$$
We set
$$
\widetilde \CM_2(M;\beta) : = \bigcup_{J \in \CJ_\omega} \widetilde \CM_2(M,J;\beta).
$$
We have the natural projection map $\pi:\widetilde \CM_2(M;\beta) \to \CJ_\omega$.

We have the natural evaluation map
$$
ev: \CF_2(\Sigma,M;\beta) \to M \times M \, ; \quad ev((j,u),(z_1,z_2)) =
(u(z_1),u(z_2)).
$$
Then the above map $\Upsilon_2$ defines a map
$$
\Upsilon_2 : \CJ_\omega \times \CF_2(\Sigma,M;\beta) \to \CH'' \times M
\times M.
$$
We now prove the following lemma by a standard argument via the linearization
of $\Upsilon_2$.

\begin{prop}\label{trans-Upsilon2} The map $\Upsilon_2$ is transverse to the
submanifold
\beastar
o_{\CH''} \times \Delta  \subset \CH'' \times (M \times M).
\eeastar
In particular the set
$$
\Upsilon_2^{-1}(o_{\CH''} \times \Delta)
$$
is a submanifold of
$$
\widetilde \CM_2(M;\beta)
$$
of codimension $2n$.
\end{prop}
\begin{proof} It is easy to check the statement on the codimension and so
we will focus on proving the submanifold property.

The linearization of $\Upsilon_2$ is given by
\be \label{eq:DUpsilon2}
(B, (\xi,v_1,v_2))
\mapsto (D_{J,u}\delbar (B,\xi), \xi(z_1) + du(v_1), \xi(z_2) +
du(v_2)).
\ee
We will focus on the problem of finite dimensional
transversality of the linear map
$$
(\xi, v_1,v_2) \mapsto (\xi(z_1) + du(v_1), \xi(z_2) + du(v_2))
$$
to the subspace $T\Delta \subset T(M\times M)$. This diagonal
transversality is well-known (e.g., see Proposition 3.4.2
\cite{mcd-sal04}).
In fact an easier variation of our proof of the 1-jet evaluation
transversality adapted to the usual 0-jet evaluation map gives rise to a
simple proof of the well-known transversality result of the evaluation map
e.g. of Proposition 3.4.2 \cite{mcd-sal04}. For the readers' convenience, we
provide the details of this transversality result in the Appendix.

Now we consider the natural projection
$$
\pi_{\Upsilon_2} : \Upsilon_2^{-1}(o_{\CH''} \times \Delta)
\to \CJ_\omega
$$
which is the restriction of the projection map
$
\pi: \widetilde\CM_2(M;\beta)\to \CJ_\omega.
$
Since $\pi$ has the index $2(c_1(\beta) + (3-n)(g-1)) + 4$,
the Fredholm index of $\pi_{\Upsilon_2}$ is given by $2(c_1(\beta) + (3-n)(g-1)) + 4 - 2n$.

Therefore for any regular value $J$ of $\pi_{\Upsilon_2}$,
$$
\widetilde\CM_2^{doub}(M,J;\b):=
\Upsilon_2^{-1}(o_{\CH^{''}}\times \Delta)\cap \pi^{-1}(J)
$$
is a smooth manifold of dimension
$$
2(c_1(\beta) +n(1-g))+ 4-2n.
$$
We just set
$$
\CJ_\omega^{inj} = \mbox{the set of regular values of $\pi_{\Upsilon_2}$}.
$$
Again we define
$$
\CM_2^{doub}(M,J;\b):=\widetilde\CM_2^{doub}(M,J;\b)/\Aut(\Sigma),
$$
where $\Aut(\S)$ acts on marked Riemann surfaces $((\S,j),(z_1,z_2))$
by conformal equivalence and then on maps from them.
Geometrically, $\CM_2^{doub}(M,J;\b)$ consists of $J$-holomorphic
curves in class $\b$ with self-intersections. It is a smooth orbifold of
dimension
$$
\text{dim}\CM_2^{doub}(M,J;\b)=2\big( c_1(\beta) +(3-n)(g-1)+2-n
\big).
$$
Therefore for any $J \in \CJ_\omega^{inj}$
$\CM_2^{doub}(M,J;\b)$ will be empty whenever
$$
c_1(\beta) + (3-n)(g-1)< n - 2
$$
and in particular when $n\ge 3$ and $c_1(\beta) \le 0$.
\end{proof}

We denote
\beastar
\CJ_\omega^{emb} & = & \CJ_\omega^{imm} \cap
\CJ_\omega^{inj}
\eeastar
which is again of second category of $\CJ_\omega$ since both $\CJ_\omega^{imm}$
and $\CJ_\omega^{inj}$ are of second category thereof.

We summarize the discussion in this
section into the following theorem

\begin{thm}\label{thm:embedding} There exists a subset $\CJ_\omega^{emb} \subset
\CJ_\omega$ of second category such that for $J \in
\CJ_\omega^{emb}$, all somewhere injective $(j,J)$-holomorphic maps
$u:\Sigma \to M$ are embeddings for any $j \in \CM_g$, provided
\be\label{eq:embedded} c_1(\beta) + (3-n)(g-1)< n - 2. \ee
\end{thm}

\begin{proof}[Proof of Theorem \ref{thm:embedded}]
Theorem \ref{thm:embedded} immediately follows from Theorem \ref{thm:immersion}
and Theorem \ref{thm:embedding} by dimension counting.
We have only to note that when $c_1(\beta) \leq 0$ and $n \geq 3$, both
inequalities \eqref{eq:immersed} and \eqref{eq:embedded} are satisfied.
\end{proof}

\section{Compactification of moduli spaces in Calabi-Yau threefolds}
\label{sec:compactification}

In this section, we restrict our attention to the case
$n = 3$, $c_1 = 0$ and let $J \in \CJ_\omega^{emb}$.
By the dimension counting argument using the evaluation maps similar
to the one in section \ref{sec:embedding}, the following is easy to
prove (See the proof of Theorem 6.3.1 \cite{mcd-sal04}
or \cite{oh-zhu} for the set-up of the proof in a somewhat
different context).

\begin{lem}\label{eq:no-nodal}
There exists a subset $\CJ_\omega^{nodal} \subset
\CJ_\omega^{emb}$ of $\CJ_\omega$ of second category such that any
two somewhere injective curves do not intersect unless they have
identical images.
\end{lem}
\begin{proof}
Recall \eqref{vdim} that the virtual dimension of $\CM(M,J;\beta)$
is given by $(n-3)(1-g)$ for $(M,\omega)$ with $c_1 = 0$ and so is
zero for $n = 3$. This result follows by the standard dimension
counting argument for Fredholm regular \emph{somewhere injective}
curves.
\end{proof}

This rules out nodal degeneration in the Gromov compactification
and gives rise to the following compactification result.

\begin{thm}\label{prop:structure} Let $J \in \CJ_\omega^{nodal}$.
Fix $\beta \in H_2(M,\Z)$ and consider a sequence of smooth
$J$-holomorphic maps $(u_i,\Sigma_i)$ with $u_i : \Sigma_i \to M$ in
class $[u_i] = \beta$. Suppose that $(u,\Sigma)$ is its stable
limit. Then there exists an integer $d\geq 1$ a class $\gamma \in
H_2(M,\Z)$ with $\beta = d \gamma$ and an embedded curve $u':
\Sigma' \to M$ with $[u'] = \gamma$ such that $(u,\Sigma)$ factors
as the composition $u = u' \circ \phi$ where $\phi : \Sigma \to
\Sigma'$ is a stable map into $\Sigma'$. Moreover, in the latter
case, each non-constant irreducible component of $\Sigma$ is a
ramified covering of $\Sigma'$.
\end{thm}
\begin{proof}  If a $J$-holomorphic curve is somewhere injective
and Fredholm regular , by dimension formula it is isolated and so
there is no $J$-holomorphic curves in a sufficiently small
$C^\infty$-neighborhood thereof. Furthermore from the above lemma,
they do not intersect unless their images coincide whenever $J \in \CJ^{nodal}_\omega$.

Combining these, we derive that every irreducible component of
$(u,\Sigma)$ either is constant or has its image coinciding with,
say, that of an embedded curve $u': \Sigma' \to M$. Denote by $C$
the image of $u'$. Then $(u,\Sigma)$ defines a stable map into the
curve $C$. Since $u': \Sigma' \to C$ is a biholomorphic map, the
stable map $u$ induces one into $\Sigma'$. Denoting this stable
map by $\phi: \Sigma \to \Sigma'$, we have finished the proof.
\end{proof}

\section{Appendix : A proof of 0-jet evaluation map transversality}

In this appendix, we give a conceptually natural proof of the following well-known
evaluation map transversality by adapting the proof of the 1-jet
evaluation transversality given in the present paper.
A detailed proof of this evaluation transversality is given in the proof of Proposition
3.4.2 \cite{mcd-sal04}, with which our proof given here would like to be compared.

\begin{thm}[0-jet evaluation transversality]
We consider the map
$$
\Upsilon_0(J,(j,u),z) = (\delbar(J,(j,u)), u(z))
$$
as a map from $\CJ_\omega \times \CF_1(\Sigma,M;\beta) \to \CH'' \times M$.
Then $\Upsilon_0$ is transverse to the submanifold
$$
o_{\CH''} \times \{p\} \subset \CH'' \times M
$$
for any given point $p \in M$.
\end{thm}
\begin{proof} Its linearization $D\Upsilon_0(J,(j,u),z)$ is given by
$$
(B,(b,\xi),v) \mapsto \left(D_{J,(j,u)}\delbar(B,(b,\xi)), \xi(u(z)) + du(z)(v)\right)
$$
for $B \in T_J\CJ_\omega, \, b \in T_j\CM(\Sigma), \, v \in T_z\Sigma$ and
$\xi \in T_u\CF(\Sigma, M;\beta)$. This defines a linear map
$$
T_J\CJ_\omega \times T_{((j,u),z)}\CF_1(\Sigma,M;\beta)
 \to \Omega_{(j,J)}^{(0,1)}(u^*TM) \times
T_{u(z)}M.
$$
We would like to prove that this linear map is surjective at every
element $(u,z_0) \in \CF_1(\Sigma,M;\beta)$ i.e., at the pair
$(u,z_0)$ satisfying
$$
\delbar_{(j,J)}u = 0, \quad u(z_0) = p.
$$
For this purpose, we need to study solvability of the system of equations
\be
D_{J,(j,u)}\delbar(B,(b,\xi)) = \gamma, \quad \xi(u(z_0)) + du(v) =
X_0 \ee
for given $\gamma \in \Omega_{(j,J)}^{(0,1)}(u^*TM)$ and $X_0
\in T_{u(z_0)}M$. Again it will be enough to consider the case $b = 0
= v$. Then this equation is reduced to
\be\label{eq:b=0=v}
D_{J,u}\delbar(B,\xi)) = \gamma, \quad \xi(u(z_0)) = X_0.
\ee
For the
map $\Upsilon_0$ to be differentiable, we need to choose the
completion of $\CF(\Sigma,M;\beta)$ in the $W^{k,p}$-norm for $k
\geq 2$.

Now we study \eqref{eq:b=0=v} for $\xi \in W^{2,p}$ similarly as in sections
\ref{sec:immersion} and \ref{appendix} : This time we can use $W^{2,p}$-norm
instead of $W^{3,p}$-norm since the 0-jet evaluation map does not involve
taking a derivative of the map $u$ unlike the 1-jet evaluation map.
We regard
$$
\Omega^{(0,1)}_{1,p}(u^*TM) \times T_{u(z_0)}M : = \CB_0
$$
as a Banach space with the norm $\|\cdot \|_{1,p} + |\cdot|$ similarly as
before with $H^{(1,0)}_{(J,(j,u),z_0)}$ replaced by $T_{u(z_0)}M$.
By the same reasoning, we apply the Fredholm alternative and study those
$(\eta, X_0)$ that satisfy the equation
$$
\langle D_u\delbar_{(j,J)}\xi + \frac{1}{2}B\circ du \circ j, \eta \rangle
+ \langle \xi, \delta_{z_0} X_0 \rangle = 0
$$
for all $B$ and $\xi$ of $C^\infty$ where $\delta_{z_0}$ is the
Dirac-delta function supported at $z_0$.

Now the rest of the proof will duplicate the proofs of Proposition
\ref{trans-Upsilon} and Lemma \ref{lem:eta=0} with $\del \xi$ replaced by $\xi$.
This finishes the proof of solvability of
\eqref{eq:b=0=v} for any given $\gamma \in W^{1,p}$ and $X_0 \in T_{u(z_0)}M$.
As before  if $\gamma \in W^{k-1,p}$, then
$\xi \in W^{k,p}$ by the elliptic regularity. This finishes the proof of surjectivity of the map
\beastar
&{}& \quad (B,(0,\xi),0) \mapsto (D_{J,u}\delbar(B,\xi)), \xi(u(z))) : \\
&{}&T_J\CJ_\omega \times \Omega^0_{k,p}(u^*TM)
\to \Omega^{(0,1)}_{k-1,p}(u^*TM) \times T_{u(z_0)}M
\eeastar
and hence proves the required transversality.
\end{proof}

A proof of the diagonal transversality, which is the transversality of the map
\eqref{eq:DUpsilon2}, can be given by an obvious modification of
the above proof by considering the map
$$
\Upsilon_2 : \CJ_\omega \times \CF_2(\Sigma,M;\beta) \to \CH'' \times M \times M
$$
whose details we omit here.


\begin{thebibliography}{HLS}

\bibitem[F]{floer} Floer, A., {\em The unregularized gradient flow of the symplectic
action,} Comm. Pure Appl. Math. 41 (1988), 775 - 813.

\bibitem[G]{gromov} Gromov, M., {\em Pseudo-holomorphic curves in
symplectic manifolds,} Invent. Math. 82 (1985), 307 - 347.

\bibitem[GS]{gelfand} Gelfand, I.M., Shilov, G.E., Generalized
Functions, vol 2, Academic Press, New York and London, 1968

\bibitem[Ho]{Ho} H\"ormander, L., The Analysis of Linera
Differential Operators II, Compre. Studies in Math. 257, Springer-Verlag,
1983, Berlin.

\bibitem[HT]{hutch-taubes} Hutchings, M., Taubes, C.H., {\em Gluing pseudo-holomorphic
curves along branched covered cylinders II}, J. Symplectic Geom. 5 (2007), 43--137,
 math.SG/0705.2074.

\bibitem[K]{konts} Kontsevich, M., A talk in the conference
``Mirror symmetry and related topics'' in University of Miami, January 2008.

\bibitem[M]{mcduff} McDuff D., {\em Example of symplectic structures,}
Invent. Math., 89 (1987), 12--36

\bibitem[MS]{mcd-sal04} McDuff D., Salamon, D., $J$-Holomorpic Curves and Symplectic Topology,
Colloquim Publications, vol 52, AMS, Providence RI, 2004.

\bibitem[OZ]{oh-zhu} Oh, Y.-G., Zhu, K., {\em Floer trajectories with immersed nodes
and scale-dependent gluing}, submitted, arXiv:0711.4187.

\bibitem[P]{pandh} Pandharipande, R., {\em Hodge integrals and degenerate contributions,}
Commun. Math. Phys. 208 (1999), 489 -- 506.

\bibitem[Si]{sikorav} Sikorav, J.-C., {\em Some properties of holomorpic
curves in almost complex manifolds}, 165 - 189, ``Holomorphic Curves in Symplectic
Geometry'', Audin, M. and Lafontaine, J. ed, Birkh\"auser, Basel, 1994

\bibitem[Sm]{smale} Smale, S., {\em
An infinite dimensional version of Sard's theorem},  Amer. J. Math. 87 (1965), 861--866

\bibitem[Wen]{wendl} Wendl, C., {\em Automatic transversality and orbifolds of
punctured holomorphic curves in dimension 4}, arXiv:0802.3842v1.


\end{thebibliography}
\end{document}